\documentclass[11pt]{article}
\usepackage{mathrsfs}
\usepackage{latexsym,lineno}
\usepackage{epsfig}
\usepackage{color}
\usepackage{amsmath}\usepackage{fleqn}\usepackage{verbatim}\usepackage{epsf}
\usepackage{amsthm}\usepackage{graphicx, float}\usepackage{graphicx}
\usepackage{amsfonts}\usepackage{amssymb}\usepackage{graphpap}
\usepackage{epic}\usepackage{curves}
\newtheorem{theorem}{Theorem}[section]

\newtheorem{corollary}{Corollary}[section]
\newtheorem{lemma}{Lemma}[section]
\newtheorem{remark}{Remark}[section]
\newtheorem{prob}{Problem}[section]

\title{\bf The distance domination of generalized de Bruijn and Kautz digraphs}
\author{Yanxia Dong$^{1}$, \, Erfang Shan$^{2}$\thanks{{\em Corresponding author.\
Email address}: efshan@shu.edu.cn(E.F. Shan)}, \, Xiao Min$^3$\\
{\small $^1$Department of Mathematics, Shanghai University, Shanghai
200444, P.R. China}\\
{\small $^2$School of Management, Shanghai University, Shanghai
200444, P.R. China}\\
{\small $^3$College of Mathematics, Physics and Information Engineering, Jiaxing University,} \\
{\small Zhejiang 314001, P.R. China}\\}

\date{}
\begin{document}
\maketitle \baselineskip 17pt
\begin{abstract}
\baselineskip 17pt
Let   $G=(V,A)$ be a digraph and $k\ge 1$ an integer. For $u,v\in V$, we say that the vertex $u$   distance $k$-dominate
 $v$ if the distance from $u$ to $v$ at most $k$. A set $D$ of vertices in $G$ is a distance $k$-dominating set if for each vertex of $V\setminus D$ is distance $k$-dominated by some vertex of $D$. The {\em distance $k$-domination number} of $G$, denoted by
$\gamma_{k}(G)$, is the minimum cardinality of a distance $k$-dominating set of $G$.
Generalized de Bruijn digraphs $G_B(n,d)$ and generalized Kautz digraphs $G_K(n,d)$ are  good candidates for interconnection networks.
Tian and Xu showed
that $\big \lceil n\big/\sum_{j=0}^kd^j\big\rceil\le
\gamma_{k}(G_B(n,d))\le \big\lceil n/d^{k}\big\rceil$ and $\big \lceil n \big/\sum_{j=0}^kd^j\big\rceil\le
\gamma_{k}(G_K(n,d))\le \big\lceil n/d^{k}\big\rceil$.
In this paper we prove that
every generalized de Bruijn digraph $G_B(n,d)$ has the distance $k$-domination number $\big\lceil n\big/\sum_{j=0}^kd^j\big\rceil$ or $\big\lceil n\big/\sum_{j=0}^kd^j\big\rceil+1$, and  the distance $k$-domination number of every generalized Kautz digraph $G_K(n,d)$
bounded above by $\big\lceil n\big/(d^{k-1}+d^{k})\big\rceil$.
Additionally, we present various sufficient conditions for $\gamma_{k}(G_B(n,d))=\big\lceil n\big/\sum_{j=0}^kd^j\big\rceil$ and $\gamma_{k}(G_K(n,d))=\big\lceil n\big/\sum_{j=0}^kd^j\big\rceil$.

\bigskip

\end{abstract}

\section{Introduction}
In this paper we deal with  directed graphs (or digraphs) which admit self-loops
but no multiple arcs. Unless otherwise defined, we follow \cite{ch, glp} for terminology and definitions.
Let $G$ be a digraph with  {\em
vertex set} $V(G)$ and  {\em arc set} $A(G)$. If there is an arc from $u$
to $v$, i.e., $(u,v)\in A(G)$, then $v$ is called an {\em out-neighbor} of $u$; we also say that $u$ {\em dominates} $v$.
%And  $u$ is called a {\em predecessor} of $v$ and $v$ is called a {\em successor} of $u$.
The {\em out-neighborhood} $O(u)$ of a vertex $u$ is the set $\{v: (u,v)\in A(G)\}$. For $S\subseteq V(G)$, its {\em out-neighborhood} $O(S)$ is the set $\cup_{u\in S}O(u)$. Set $O_0(u)=\{u\}$ and $O_{1}(u)=O(u)$,
we define recursively  $O_{i}(u)$, called {\em $i$-th out-neighborhood} of $u$,  by $\{O(O_{i-1}(u))\}$ for $i\geq1$.
The {\em $i$-th out-neighborhood} of $S$ is  the set  $O_{i}(S)=\cup_{u\in S}O_{i}(u)$. The {\em closed out-neighborhood} $O[u]$ of $u$ is the set $O(u)\cup\{u\}$, and $O[S]$ and $O_{i}[S]$ are defined analogously.

For $x, y\in V(G)$, the {\em distance} $d_{G}(x,y)$ from $x$ to  $y$ is the length of an shortest $(x,y)$-directed path in $G$.
%The maximum distance between any two vertices in $G$ is called the diameter of $G$.
Let $k$ be a positive integer.
A subset $D\subseteq V(G)$ is called a {\em distance $k$-dominating set} of $G$ if for every vertex $v$ of $V(G)\setminus D$, there is a vertex $u\in D$ such that $d_G(u,v)\le k$, i.e., $\cup_{i=0}^{k}O_i(D)=V(G)$.
 The {\em distance $k$-domination number} of $G$,
denoted by $\gamma_{k}(G)$, is the minimum cardinality of a distance $k$-dominating set of $G$.
In particular, the distance 1-dominating set  is the ordinary dominating set, which has been well studied \cite{hhs}.

%From the definition, $\gamma_{k}(G)=1$ for $k$ not less than the diameter of $G$. Thus we assume $k$ be any positive %integer less than the diameter of $G$.
%For a vertex $u\in V(G)$, it is clearly that $d_{G}(u,v)\le k$ for $v\in O_{k}(u)$.

Slater \cite{hhs} termed a distance $k$-dominating set as a $k$-basis and also gave an interpretation for a $k$-basis in terms of communication networks. Since then many researchers pay much attention to this subject, for example \cite{fv,s,sse}. The concept of distance domination in graphs finds applications in many structures and situations which give rise to graphs. A minimum distance $k$-dominating set of $G$ may be used locate a minimum number of facilities (such as utilities, police stations, hospitals, transmission towers, blood banks, waste disposal dump) such that every intersection is within $k$ city block of a facility. Barkauskas and Host \cite{a} showed that the problem of determining $\gamma(G)$ is NP-hard for a general graph.

The network topology has a great impact on the system performance and reliability \cite{xu}.
There are some well-known networks with good properties such as de Bruijn networks, Kautz networks and their generalizations
(see, for example, \cite{bp, det, dew, hkmp, xu}).
Generalized de Bruijn and Kautz networks, denoted by $G_B(n,d)$ and $G_K(n,d)$ respectively, were introduced by Imase and Itoh \cite{ii}. The generalization
removes the restriction on the cardinality of vertex set and make the network more general and valuable as a network model.
A lot of features make it suitable for implementation of reliable networks. The most
important feature such as small diameter \cite{ii}, high connectivity \cite{iso}, easy routing, and high reliability.

The generalized de Bruijn digraph $G_{B}(n,d)$ is defined by
congruence equations as follows:
$$
\left\{\begin{array}{lll}
 V(G_{B}(n,d))= \{0,1,2,\ldots,n-1 \}
 \\
 A(G_{B}(n,d))=\{(x,y)\!\mid\! y\equiv dx+i \,($mod$ ~n), 0\leq i\leq d-1
 \}.
\end{array}\right.
$$
In particular, if $n=d^{m}$, then $G_{B}(n,d)$ is the de
Bruijn digraph $B(d,m)$. The generalized Kautz digraph $G_{K}(n,d)$
is defined by following congruence equation:
$$
\left\{\begin{array}{lll}
 V(G_{K}(n,d))= \{0,1,2,\ldots,n-1 \}
 \\A(G_{K}(n,d))=\{(x,y)\!\mid\! y\equiv -dx-i \,($mod$ ~n),1\leq i\leq d
 \}.
\end{array}\right.
$$

In
particular, if $n=d^{m}+ d^{m-1}$, then $G_{K}(n,d)$ is the Kautz
digraph $K(d,m)$. \noindent The graphs $G_{B}(6,3)$ and $G_{K}(9,2)$
are exhibited in Fig. 1.

\begin{figure}[h]
\setlength{\unitlength}{.015in}
\begin{picture}(360,160)
\put(30,5)

\multiput(50,60)(40,0){3}{\circle*{5}}
\multiput(50,140)(40,0){3}{\circle*{5}}

\multiput(50,140)(40,0){2}{\line(1,0){40}}
\multiput(90,60)(40,0){2}{\line(-1,0){40}}
\multiput(50,140)(40,0){2}{\vector(1,0){20}}
\multiput(90,60)(40,0){2}{\vector(-1,0){20}}

\put(90,60){\vector(-1,2){30}} \put(50,60){\vector(1,2){30}}
\put(90,140){\vector(1,-2){30}} \put(130,140){\vector(-1,-2){30}}
\put(90,60){\line(-1,2){40}} \put(50,60){\line(1,2){40}}
\put(90,140){\line(1,-2){40}} \put(130,140){\line(-1,-2){40}}

\qbezier(50,60)(40,100)(50,140)\qbezier(50,60)(60,100)(50,140)
\qbezier(90,60)(80,100)(90,140)\qbezier(90,60)(100,100)(90,140)
\qbezier(130,60)(120,100)(130,140)\qbezier(130,60)(140,100)(130,140)

\qbezier(50,60)(40,50)(50,40)\qbezier(50,60)(60,50)(50,40)
\qbezier(130,60)(120,50)(130,40)\qbezier(130,60)(140,50)(130,40)
\qbezier(50,140)(40,150)(50,160)\qbezier(50,140)(60,150)(50,160)
\qbezier(130,140)(120,150)(130,160)\qbezier(130,140)(140,150)(130,160)

\multiput(45,100)(40,0){3}{\vector(0,1){0.5}}
\multiput(55,100)(40,0){3}{\vector(0,-1){0.5}}
\multiput(55,50)(80,0){2}{\vector(0,1){0.5}}
\multiput(55,150)(80,0){2}{\vector(0,-1){0.5}}

\put(38,135){0}\put(90,148){1}\put(140,135){3}
\put(38,55){2}\put(90,47){4}\put(140,55){5}

\put(50,5){Figure 1 (a):  $G_{B}(6, 3)$}

\multiput(240,50)(40,0){3}{\circle*{5}}
\multiput(240,150)(40,0){3}{\circle*{5}}
\multiput(300,90)(0,20){2}{\circle*{5}} \put(360,90){\circle*{5}}

\multiput(240,150)(40,0){2}{\line(1,0){40}}
\multiput(240,50)(40,0){2}{\line(1,0){40}}
\multiput(240,150)(40,0){2}{\vector(1,0){20}}
\multiput(240,50)(40,0){2}{\vector(1,0){20}}

\put(280,150){\line(1,-2){20}} \put(280,50){\line(1,2){20}}
\put(300,110){\line(-1,-3){20}} \put(300,90){\line(-1,3){20}}
\put(320,150){\line(2,-3){40}} \put(320,50){\line(1,1){40}}
\put(280,150){\vector(1,-2){10}} \put(280,50){\vector(1,2){10}}
\put(300,110){\vector(-1,-3){10}} \put(300,90){\vector(-1,3){10}}
\put(320,150){\vector(2,-3){20}} \put(320,50){\vector(1,1){20}}

\qbezier(240,50)(220,100)(240,150)\qbezier(240,50)(260,100)(240,150)
\qbezier(300,110)(295,130)(320,150)\qbezier(320,150)(325,130)(300,110)
\qbezier(300,90)(295,70)(320,50)\qbezier(320,50)(325,70)(300,90)
\put(230,100){\vector(0,1){0.5}}\put(250,100){\vector(0,-1){0.5}}
\put(305,135){\vector(1,1){0.5}}\put(315,125){\vector(-1,-1){0.5}}
\put(303,70){\vector(1,-1){0.5}}\put(317,70){\vector(-1,1){0.5}}

\qbezier(360,90)(355,190)(240,150)\qbezier(360,90)(350,10)(240,50)
\put(342,144.5){\vector(-1,1){0.5}}\put(340.5,49){\vector(-1,-1){0.5}}

\put(230,45){8}\put(280,40){1}\put(320,40){6}
\put(230,145){0}\put(282,151){7}\put(325,145){2}
\put(305,105){3}\put(305,87){5}\put(365,85){4}

\put(250,5){Figure 1 (b):  $G_{K}(9, 2)$}

\end{picture}

\end{figure}

The structure properties of the generalized de Bruijn  and  Kautz digraphs  receive more attention.
Du et al. \cite{dhhz} studied the hamiltonian property of generalized de Bruijn and Kautz networks.
Also, several structural objects such as spanning trees, Eulerian tours \cite{lz}, closed walks \cite{sso} and small
cycles \cite{hs} have been counted.
Shan et al. \cite{sck,sd,sdc} studied the absorbants and twin domination of generalized de Bruijn digraphs. Recently, Dong  et al. \cite{dsk} completely determined  the domination number of generalized de Bruijn digraphs. Wang \cite{wa}
showed that there is an  efficient twin dominating set in $G_B(n,d)$ with $n=c(d+1)$ if and only if $d$ is even and relatively prime to $c$.
More studied progress on the generalized de Bruijn and Kautz networks can be found in \cite{duc,tx,xu}.

In order to make our arguments easier to follow we introduce the {\em modulo interval} so as to
represent the out-neighborhood of each vertex in $G_{B}(n,d)$ and $G_{K}(n,d)$. Let $I=\{0, 1, \cdots, n-1\}$
denote the vertex set of $G_{B}(n,d)$.  For any integers $i,j$ satisfying $i\not \equiv j$ (mod $n$), a {\em modulo interval} $[i,j]$ $({\rm mod} ~n)$, with respect to modulo $n$, is defined by
\begin{equation*}
[i,j] \, ({\rm mod} ~n)=
\begin{cases}
\{i,i+1, \ldots, j\}\, ({\rm mod} ~n) & \text{if $i$ (mod $n)<j$ (mod $n$)},\\
\{i, \ldots, n-1, 0, \ldots, j\} \,({\rm mod} ~n) &  \text{if $i$ (mod $n)>j$ (mod $n$)}.
\end{cases}
\end{equation*}
By the definitions,  $I=[0, n-1]$, and for each $j\in [0,n-1]$, clearly $O(j)=[jd,jd+(d-1)]$ $({\rm mod} ~n)$ in $G_{B}(n,d)$ and $O(j)=[-jd-d,-jd-1]$ $({\rm mod} ~n)$ in $G_{K}(n,d)$.

Notice that if $d=1$ then the graph $G_B(n,1)$ (or $G_K(n,1)$) has $n$ self-loops.  Throughout this paper, we  always assume $d\ge 2$ and $n\ge d$. If the set $D=\{x,x + 1, \cdots, x + k\}$ ({\rm mod} $n$) is a dominating set or a distance $k$-dominating set of
$G_{B}(n,d)$ (or $G_K(n,d)$), then $D$ is called a {\em consecutive dominating set} or a {\em consecutive distance $k$-dominating set} of $G_{B}(n,d)$ (or $G_K(n,d)$).
 A {\em  consecutive minimum dominating set} of $G_{B}(n,d)$ (or $G_K(n,d)$) is a consecutive dominating set with cardinality  $\gamma(G_{B}(n,d))$ (or $\gamma(G_{K}(n,d))$) and a {\em consecutive distance $k$-dominating set} of $G_{B}(n,d)$ (or $G_{K}(n,d)$) is a consecutive distance $k$-dominating set with cardinality  $\gamma_{k}(G_{B}(n,d))$ (or $\gamma_{k}(G_{K}(n,d))$).

Tian and Xu \cite{tx} established the upper and lower bounds on the distance $k$-domination number of $G_{B}(n,d)$ and $G_{K}(n,d)$. This paper continues to study distance $k$-domination in generalized de
Bruijn and Kautz digraphs. In
Subsection 2.1, we show that  every generalized de Bruijn digraph $G_{B}(n,d)$ has the distance $k$-domination number either $\big\lceil n\big/\sum_{j=0}^kd^{j}\big\rceil$ or $\big\lceil n\big/\sum_{j=0}^kd^{j}\big\rceil+1$. In Subsection 2.2, we derive various sufficient conditions for $\gamma_{k}(G_B(n,d))=\big\lceil n\big/\sum_{j=0}^kd^{j}\big\rceil$.
In Section 3, we gives a sharp upper bound of $\gamma_{k}(G_{K}(n,d))$, which improves the previous upper bound of $\gamma_{k}(G_{K}(n,d))$, due to Tian and Xu \cite{tx}.  In closing section, we pose two open problems.

\section{The minimum distance $k$-dominating sets in $G_B(n,d)$}
In the first subsection of this section, by constructing a distance $k$-dominating set of an arbitrary generalized de Bruijn digraph $G_B(n,d)$, we  show that the distance $k$-domination number of $G_B(n,d)$ has exactly two values.
In next subsection, we describe various sufficient conditions for the distance $k$-domination number
equal to one of two values.

\subsection{The distance $k$-domination number of $G_B(n,d)$}

Tian and Xu \cite{tx} observed the following upper and lower bounds on $\gamma_{k}(G_B(n,d))$.

\begin{lemma}{\rm(\cite{tx})}\label{lem2.2}
 For every generalized de Bruijn digraph $G_B(n,d)$,
 $$\Big\lceil n\Big/ \sum_{j=0}^kd^{j}\Big\rceil\le
\gamma_{k}(G_B(n,d))\le \Big\lceil\frac{n}{d^{k}}\Big\rceil.$$
\end{lemma}

We are ready to improve the above upper bound on $\gamma_{k}(G_B(n,d))$ by directly constructing a (consecutive) distance $k$-dominating set of $G_B(n,d)$ with cardinality $\big\lceil n/(\sum_{j=0}^kd^j)\big\rceil+1$.
The following lemma  plays a key
role in constructing such a distance $k$-dominating set of $G_B(n,d)$.

\begin{lemma} \label{lem2.20}
 Every generalized de Bruijn digraph $G_B(n,d)$ contains a vertex $x$ satisfying the following inequality:
 \begin{eqnarray}\label{formu1}
x+\Big\lceil n\Big/\sum_{j=0}^kd^{j}\Big\rceil-(d-2)\leq dx\leq x+\Big\lceil n\Big/\sum_{j=0}^kd^{j}\Big\rceil\, \,({\rm mod}\, n).
\end{eqnarray}
\end{lemma}
\noindent {\bf Proof.}  We choose an arbitrary vertex $x\in V(G_B(n,d))$. If $x$ satisfies (\ref{formu1}), we are done. Otherwise, the vertex $x$ clearly satisfies either
\begin{align*}
0\le dx\le x+\Big\lceil n\big/\sum_{j=0}^kd^j\Big\rceil-(d-1)\,  \,({\rm mod}\, n)
\end{align*}
or
\begin{align*}
x+\Big\lceil n\big/\sum_{j=0}^kd^j\Big\rceil+1\le dx\le n-1 \, \,({\rm mod}\, n).
\end{align*}

We find the desired vertex by distinguishing the following two cases.

{\em Case} 1. $0\le dx\le x+\big\lceil n\big/\sum_{j=0}^kd^j\big\rceil-(d-1)\, \,({\rm mod}\, n) $. Note that if $x$ increases by integer $i$, then the value of $dx$ is increased  to $d(x+i)=dx+di$. In this case, we find the
desired vertex  by increasing the value of $x$. Since $dx\le x+\big\lceil n\big/\sum_{j=0}^kd^j\big\rceil-(d-1)\, \,({\rm mod}\, n)$, there exists an  integer
$i$ ($\ge 0$) such that $x$ and $i$ satisfy the following inequality
\begin{eqnarray}\label{formu2}
 d(x+i)\le x+\Big\lceil n\Big/\sum_{j=0}^kd^j\Big\rceil-(d-2) \, \,({\rm mod}\, n),
\end{eqnarray}
since $i=0$ satisfies the inequality.
  Let $i$ be the maximal integer  satisfying (\ref{formu2}).
We claim that
\begin{eqnarray}\label{formu3}
d(x+i)\ge (x+i)+\Big\lceil n\Big/\sum_{j=0}^kd^j\Big\rceil-2(d-2) \, \,({\rm mod}\, n).
\end{eqnarray}
Indeed, if $d(x+i)\le (x+i)+\left\lceil n\big/\sum_{j=0}^kd^j\right\rceil-2(d-2)-1\, \,({\rm mod}\, n)$, then
\begin{eqnarray*}
d(x+i+1)\le (x+i+1)+\Big\lceil n\Big/\sum_{j=0}^kd^j\Big\rceil-(d-2)\, \,({\rm mod}\, n).
\end{eqnarray*}
So $i+1$ satisfies (\ref{formu2}) too, this contradicts the maximality of $i$. Hence (\ref{formu3}) follows.
 If the equality holds in (\ref{formu2}), that is,
\begin{eqnarray*}
d(x+i)=x+\Big\lceil n\Big/\sum_{j=0}^kd^j\Big\rceil-(d-2) \, \,({\rm mod}\, n),
\end{eqnarray*}
 then $x+i$ satisfies (\ref{formu1}).
So we replace $x$ by $x+i$, and obtain the desired vertex. Otherwise, by (\ref{formu3}), we have
\begin{eqnarray*}
(x+i)+\Big\lceil n\Big/\sum_{j=0}^kd^j\Big\rceil-2(d-2)\le d(x+i)\le (x+i)+\Big\lceil n\Big/\sum_{j=0}^kd^j\Big\rceil-(d-1)\, \,({\rm mod}\, n).
\end{eqnarray*}
Hence,
\begin{eqnarray*}
(x+i+1)+\Big\lceil n\Big/\sum_{j=0}^kd^j\Big\rceil-(d-3)
\le d(x+i+1)
\le (x+i+1)+\Big\lceil n\Big/\sum_{j=0}^kd^j\Big\rceil \, \,({\rm mod}\, n).
\end{eqnarray*}
Clearly, $x+i+1$ satisfies (\ref{formu1}).  Thus we replace $x$ by $x+i+1$ and obtain the desired vertex.

{\em Case} 2. $x+\big\lceil n\big/\sum_{j=0}^kd^j\big\rceil+1\le dx\le n-1 \, \,({\rm mod}\, n)$. We can obtain the
desired vertex  by decreasing the value of $x$. Clearly, there exists an  integer
$i$ ($\ge 0$) such that $x$ and $i$ satisfy the following inequality
\begin{eqnarray}\label{formu4}
 d(x-i)\ge (x-i)+\Big\lceil n\Big/\sum_{j=0}^kd^j\Big\rceil\, \,({\rm mod}\, n),
\end{eqnarray}
since the inequality $dx\ge x+\big\lceil n\big/\sum_{j=0}^kd^j\big\rceil+1$ implies that $i=0$ satisfies (\ref{formu4}). Let $i$ be the maximal integer  satisfying (\ref{formu4}). We claim that
\begin{eqnarray}\label{formu5}
 d(x-i)\le (x-i)+\Big\lceil n\Big/\sum_{j=0}^kd^j\Big\rceil+d-2\, \,({\rm mod}\, n).
\end{eqnarray}
Suppose, to the contrary, that $d(x-i)\ge (x-i)+\big\lceil n\big/\sum_{j=0}^kd^j\big\rceil+d-1$ (mod $n$). Equivalently,
\begin{eqnarray*}
d(x-(i+1))\ge (x-(i+1))+\Big\lceil n\Big/\sum_{j=0}^kd^j\Big\rceil \,({\rm mod}\, n).
\end{eqnarray*}
But then $i+1$ satisfies (\ref{formu4}). This contradicts the maximality of $i$. Thus (\ref{formu5}) holds.
If the equality holds in (\ref{formu4}), then
the vertex $x-i$ satisfies (\ref{formu1}). So we obtain the desired vertex by replacing $x$ by $x-i$. Otherwise,
by (\ref{formu5}), we have
\begin{eqnarray*}
(x-i)+\Big\lceil n\Big/\sum_{j=0}^kd^j\Big\rceil+1\le d(x-i)\le (x-i)+\Big\lceil n\Big/\sum_{j=0}^kd^j\Big\rceil+d-2\,\,({\rm mod}\, n).
\end{eqnarray*}
Hence,
\begin{align*}
(x-(i+1))+\Big\lceil n\Big/\sum_{j=0}^kd^j\Big\rceil-(d-2)&\le d(x-(i+1))\\
&\le (x-(i+1))+\Big\lceil n\Big/\sum_{j=0}^kd^j\Big\rceil-1\,\, ({\rm mod}\, n).
\end{align*}
Hence $x-(i+1)$ satisfies (\ref{formu1}). We obtain the desired vertex by replacing $x$ by $x-(i+1)$.
~$\Box$

\begin{theorem}\label{thm2.1}
For every generalized de Bruijn digraph $G_B(n,d)$,
$$\gamma_{k}(G_B(n,d))=\Big\lceil n\Big/\sum_{j=0}^kd^{j}\Big\rceil\, \, \mbox{or} \, \, \Big\lceil n\Big/\sum_{j=0}^kd^{j}\Big\rceil+1.$$
\end{theorem}
\noindent {\bf Proof.}
By Lemma \ref{lem2.2}, it suffices to show that $\gamma(G_B(n,d))\le \big\lceil n\big/\sum_{j=0}^kd^{j}\big\rceil+1$.
The proof is by directly constructing a (consecutive) distance $k$-dominating set of $G_B(n,d)$
with cardinality $\big\lceil n/(\sum_{j=0}^kd^j)\big\rceil+1$.
By Lemma \ref{lem2.20}, there is a vertex $x$ in $G_B(n,d)$ that satisfies (\ref{formu1}).
Let $D=\big\{x,x+1,\ldots,  x+\big\lceil n\big/\sum_{j=0}^kd^j\big\rceil\big\}$. We show that $D$ is
a distance $k$-dominating set of $G_B(n,d)$. By the definition, we need to prove that $\bigcup _{i=0}^{k}O_{i}(D)=V(G_B(n,d))$.

First, we show that  the vertices of $O_{i-1}\cup O_{i}(D)$ are consecutive for all $i$, $1\le i\le k$.  The out-neighborhoods of vertices in $D$ are given as follows.
\begin{align*}
 &O(x)=\{dx,dx+1,\ldots,  dx+d-1\} \,\, ({\rm mod}\, n),\\
&O(x+1)=\{d(x+1),d(x+1)+1,\ldots, d(x+1)+d-1\}  \,\, ({\rm mod}\, n),\\
&~~~~\vdots\\
&O\Big(x+\Big \lceil n\Big/\sum_{j=0}^kd^j\Big\rceil\Big)=\Big\{d\Big(x+\Big
\lceil n\Big/\sum_{j=0}^kd^j\Big\rceil\Big),
\ldots,  d\Big(x+\Big
\lceil n\Big/\sum_{j=0}^kd^j\Big\rceil\Big)+d-1\Big\}\,\, ({\rm mod}\, n).
\end{align*}
Then $O(D)=\big[dx,d\big(x+
\big\lceil n\big/\sum_{j=0}^kd^j\big\rceil\big)+d-1\big]\,\, ({\rm mod}\, n)$. Similarly, the $i$-th out-neighborhoods $O_{i}(D)=\big[d^{i}x,d^{i}\big(x+
\big\lceil n\big/\sum_{j=0}^kd^j\big\rceil\big)+(d-1)\sum_{j=0}^i d^j\big]\,\, ({\rm mod}\, n)$ for each $i, 1\leq i\le k$. Since $x$ satisfying the inequality (\ref{formu1}),
there exists an integer $h$, $0\leq h\leq d-2$, such that $dx= x+\big\lceil n\big/\sum_{j=0}^kd^j\big\rceil-h \,\,({\rm mod}\, n)$, so we have
\begin{align*}
&d^{2}x= d\Big(x+\Big
\lceil n\Big/\sum_{j=0}^kd^j\Big\rceil\Big)-dh \,\, ({\rm mod}\, n),\\
&d^{3}x= d^{2}\Big(x+\Big
\lceil n\Big/\sum_{j=0}^kd^j\Big\rceil\Big)-d^{2}h \,\, ({\rm mod}\, n),\\
&~~~~~~~\vdots\\
&d^{k}x= d^{k-1}\Big(x+\Big
\lceil n\Big/\sum_{j=0}^kd^j\Big\rceil\Big)-d^{k-1}h\,\, ({\rm mod}\, n).
\end{align*}
Thus  $O_{i-1}(D)\cap O_{i}(D)\neq \emptyset$ for all $i, 1\leq i\le k$. This implies that
the vertices of $O_{i-1}(D)\cup O_{i}(D)$ are consecutive, since the  vertices of $O_i(D)$ are consecutive
for each $i$, $0\le i\le k$. Therefore, the vertices of  $\bigcup_{i=0}^kO_{i}(D)$ are   consecutive.

Next we show that $\bigcup _{i=0}^{k}O_{i}(D)$ contains all the vertices of $G_B(n,d)$.
Note that $O_1(D)\cap D\neq \emptyset$. Thus it suffices to
show that  $O_k(D)\cap D\neq \emptyset$.
For the last vertex in $O_k(D)$, since $x$ satisfies (1), we have
\begin{align*}
d^{k}&\Big(x+\Big\lceil n\Big/\sum_{j=0}^kd^j\Big\rceil\Big)+\big(d-1\big)\sum_{j=0}^kd^j\\
&= d^{k-1}\Big(x+\Big\lceil n\Big/\sum_{j=0}^kd^j\Big\rceil-h\Big)+d^{k}\Big
\lceil n\Big/\sum_{j=0}^kd^j\Big\rceil+\big(d-1\big)\sum_{j=0}^kd^j\\
&=d^{k-1}x+\big(d^{k}+d^{k-1}\big)\Big
\lceil n\Big/\sum_{j=0}^kd^j\Big\rceil+(d-1)d^{k}-hd^{k-1}+\big(d-1\big)\sum_{j=0}^kd^j\\
&~~~~~~\vdots\\
%\end{align*}
%\begin{align*}
&= x+\Big
\lceil n\Big/\sum_{j=0}^kd^j\Big\rceil\sum_{j=0}^kd^j -h\sum_{j=0}^{k-1}d^j+\big(d-1\big)\sum_{j=0}^kd^j\\
&= x+(d-1)+\Big
\lceil n\Big/\sum_{j=0}^kd^j\Big\rceil\sum_{j=0}^kd^j+\big(d(d-1)-h\big)\sum_{j=0}^{k-1}d^j\\
&\geq x\,\, ({\rm mod}\, n)
\end{align*}
The last inequality holds, since $d\geq2$ and
$0\leq h\leq d-2$. Hence $O_k(D)\cap D\neq \emptyset$, and so
$$\bigcup _{i=1}^{k}O_{i}(D)\supseteq \big\{x+\big
\lceil n\big/\sum_{j=0}^kd^j\big\rceil,\ldots,  n-1,0,1,\ldots,  x\big\}.$$ This implies that $\bigcup _{i=0}^{k}O_{i}(D)=V(G_B(n,d))$, that is, $D$   is a
(consecutive) distance $k$-dominating set of $G_B(n,d)$. Consequently, $\gamma_k(G_B(n,d))\le |D|=\big\lceil n\big/\sum_{j=0}^kd^j\big\rceil+1.$ ~$\Box$

For distance $k=1$ we obtain the following result.
\begin{corollary}{\rm(\cite{dsk})}\label{cor2.1}
 For every generalized de Bruijn digraph $G_B(n,d)$, either $\gamma(G_B(n,d))=\big\lceil \frac{n}{d+1}\big\rceil$ or $\gamma(G_B(n,d))=\big\lceil \frac{n}{d+1}\big\rceil+1$.
\end{corollary}

\subsection{The generalized de Bruijn digraphs $G_B(n,d)$ with $\gamma(G_B(n,d))=\big\lceil \frac{n}{d+1}\big\rceil$}
In the next subsection, we derive various sufficient conditions for the distance $k$-domination
number to achieve the value $\big\lceil n\big/\sum_{j=0}^kd^j\big\rceil$ in a generalized de Bruijn digraph $G_B(n,d)$.
\begin{theorem}\label{thm2.2}
  If there exists a vertex $x\in V(G_B(n,d))$ satisfying the following congruence equation:
\begin{eqnarray}\label{formu6}
(d-1)x\equiv \Big\lceil n\Big/\sum_{j=0}^kd^j\Big\rceil-h\, \,({\rm mod}\, n),
\end{eqnarray}
 for some $h$ where $0\leq \big(\sum_{j=0}^{k-1}d^j\big)h\leq \big(\sum_{j=0}^kd^j\big)\lceil n\big/\sum_{j=0}^kd^j\rceil-n$,
then $\gamma_{k}(G_B(n,d))=\big\lceil n\big/\sum_{j=0}^kd^j\big\rceil$, and $D=\big\{x,x+1,x+2,\ldots, x+\big\lceil n\big/\sum_{j=0}^kd^j\big\rceil-1\big\}$ is a consecutive minimum distance $k$-dominating set of $G_B(n,d)$.
\end{theorem}
\noindent {\bf Proof.}
Let $x$ be a vertex of $G_B(n,d)$ satisfying Eq. (\ref{formu6}). Note that $|D|=\big\lceil n\big/\sum_{j=0}^kd^j\big\rceil$. By Theorem \ref{thm2.1},
it is sufficient to show that $D=\{x,x+1,x+2,\ldots, x+\big\lceil n\big/\sum_{j=0}^kd^j\big\rceil-1\}$ is a distance $k$-dominating set of $G_B(n,d)$. For this purpose, we show that $\bigcup _{i=1}^{k}O_{i}(D)= V(G_B(n,d))$.

We first prove that  the vertices of $O_{i-1}(D)\cup O_{i}(D)$ are consecutive for all $i, 1\leq i\le k$.
By the definition of $G_B(n,d)$,
the out-neighborhoods $O(D)$ of $D$ are given as follows.
\begin{align*}
 &O(x)=\{dx,dx+1,\ldots,  dx+d-1\} \,\, ({\rm mod}\, n),\\
&O(x+1)=\{d(x+1),d(x+1)+1,\ldots, d(x+1)+d-1\}  \,\, ({\rm mod}\, n),\\
&~~~~\vdots\\
&O\Big(x+\Big \lceil n \Big/\sum_{j=0}^kd^j\Big\rceil-1\Big)=\Big\{d\Big(x+\Big
\lceil n \Big/\sum_{j=0}^kd^j\Big\rceil\Big)-d,d\Big(x+\Big
\lceil n \Big/\sum_{j=0}^kd^j\Big\rceil\Big)\\
&~~~~~~~~~~~~~~~~~~~~~~~~~~~~~~~~~~~~~~~~~-d+1,
\ldots,  d\Big(x+\Big
\lceil n \Big/\sum_{j=0}^kd^j\Big\rceil\Big)-1\Big\}\,\, ({\rm mod}\, n).
\end{align*}
Then $O(D)=\big[dx,dx+d\big
\lceil n \big/\sum_{j=0}^kd^j\big\rceil-1\big]\,\, ({\rm mod}\, n)$.
Similarly, we have
$O_{i}(D)=[d^{i}x, d^{i}\big(x+\big
\lceil n \big/\sum_{j=0}^kd^j\big\rceil\big)-1]\,\, ({\rm mod}\, n)$.
Clearly,  $|O_{i}(D)|=d^{i}\big
\lceil n \big/\sum_{j=0}^kd^j\big\rceil$ for all $i, 0\le i\le k$.
Since  $x$ satisfies Eq. (\ref{formu6}), we have
\begin{align*}
 O(D)&=\Big[x+\Big\lceil n \big/\sum_{j=0}^kd^j\Big\rceil-h,
d\Big(x+\Big
\lceil n \Big/\sum_{j=0}^kd^j\Big\rceil\Big)-1\Big]\,\, ({\rm mod}\, n),\\
O_{2}(D)&=\Big[d\Big(x+\Big
\lceil n \Big/\sum_{j=0}^kd^j\Big\rceil\Big)-dh,
d^{2}\Big(x+\Big
\lceil n \Big/\sum_{j=0}^kd^j\Big\rceil\Big)-1\Big] \,\, ({\rm mod}\, n),\\
&~~~~\vdots\\
O_{k}(D)&=\Big[d^{k-1}\Big(x+\Big
\lceil n \Big/\sum_{j=0}^kd^j\Big\rceil\Big)-d^{k-1}h,
d^{k}\Big(x+\Big
\lceil n \Big/\sum_{j=0}^kd^j\Big\rceil\Big)-1\Big]\,\, ({\rm mod}\, n).
\end{align*}
 Hence it can be seen that  $|O_{i-1}(D)\cap O_{i}(D)|=d^{i-1}h$ for all $i, 1\le i\le k$. Note that the vertices of each $O_{i}(D)$ ($i\ge 0$)   are consecutive.
  By the above observations, if $h=0$, then the last vertex of $O_{i-1}(D)$ and the first vertex of $O_{i}(D)$ are consecutive; while if $h>0$, then  $O_{i-1}(D)\cap O_{i}(D)\neq \emptyset$. Thus the vertices of $O_{i-1}(D)\cup O_{i}(D)$ are consecutive for all $i, 1\leq i\le k$.

We next show that $\bigcup _{i=0}^{k}O_{i}(D)=V(G_B(n,d))$.
As observed above,   we see that the vertices of $\bigcup _{i=0}^{k}O_{i}(D)$ are consecutive,
In particular, the vertices of $D\cup O_1(D)$ are consecutive.
Thus it suffices to
show that  the vertices $O_k(D)\cup D$ are consecutive.
For the last vertex in $O_k(D)$,
because $0\leq \big(\sum_{j=0}^{k-1}d^j\big)h\leq \big(\sum_{j=0}^kd^j\big)\lceil n\big/\sum_{j=0}^kd^j\rceil-n$,
we have
\begin{align*}
&d^{k}\Big(x+\Big\lceil n \Big/\sum_{j=0}^kd^j\Big\rceil\Big)-1\,\, ({\rm mod}\, n)\\
&=x+\Big(\sum_{j=0}^kd^j\Big)\Big
\lceil n \Big/\sum_{j=0}^kd^j\Big\rceil-\Big(\sum_{j=0}^{k-1}d^j\Big)h-1 \,\, ({\rm mod}\, n)\, \, (\mbox{by (\ref{formu6})})\\
&\geq x-1 \,({\rm mod}\, n).
\end{align*}
This implies that  the vertices of $O_k(D)\cup D$ are consecutive, so
$$\bigcup _{i=1}^{k}O_{i}(D)\supseteq \big\{x+\big
\lceil n\big/\sum_{j=0}^kd^j\big\rceil,\ldots,  n-1,0,1,\ldots,  x-1\big\}.$$
This implies that $\bigcup _{i=0}^{k}O_{i}(D)=V(G_B(n,d))$,
hence $D$ is a  distance $k$-dominating set of $G_B(n,d)$. This complete the proof of Theorem \ref{thm2.2}. ~$\Box$

As a special case of Theorem \ref{thm2.2}, we immediately have the following corollary.

\begin{corollary}\label{cor2.2}
Let $\sum_{j=0}^kd^j\!\mid\! n$. If there is a vertex $x\in V(G_B(n,d))$ satisfying congruence equation:
\begin{eqnarray}\label{formu7}
(d-1)x\equiv n\Big/\sum_{j=0}^kd^j\, \,({\rm mod}\, n),
\end{eqnarray}
then  $\gamma_{k}(G_B(n,d))= n\big/\sum_{j=0}^kd^j$ and $D=\{x,x+1,\cdots,x+n\big/\sum_{j=0}^kd^j-1\}$  is a consecutive minimum  distance $k$-dominating set of $G_B(n,d)$.
\end{corollary}

 \begin{remark}{\rm If $G_B(n,d)$ contains no vertex $x$ satisfying (\ref{formu6}) in Theorem \ref{thm2.2}, it is possible to encounter $\gamma_{k}(G_B(n,d))=\big\lceil n\big/\sum_{j=0}^kd^j\big\rceil+1$. For example, let $G_B(40,3)$ and $k=3$. The congruence equation $(d-1)x\equiv \big\lceil n\big/\sum_{j=0}^kd^j\big\rceil-h\, \,({\rm mod}\, n)$ is $2x\equiv 1\, \,({\rm mod}\, 40)$  where $h=0$, since  $ 40\big/\sum_{j=0}^33^j=1$. Clearly, there is no vertex satisfying $2x\equiv 1\, \,({\rm mod}\, 40)$.
We can deduce that $\gamma_{3}(G_B(40,3))=\big\lceil 40\big/\sum_{j=0}^33^j\big\rceil+1=2$. Indeed,  for each $x$ of $G_B(40,3)$, it can be verify that $\{x\}$ is not a distance $3$-dominating set of $G_B(40,3)$ by  simply enumeration.}
\end{remark}

Recalling that $G_{B}(d^{m},d)=B(d,m)$ when $n=d^m$. For  cases $k=1$ and $k=2$, the
distance $k$-domination numbers of a de Bruijn digraph $B(d,m)$ were proved by Araki \cite{a} and  Tian \cite{tx}, respectively.  As an application of Theorem \ref{thm2.2}, we can determine the
distance $k$-domination number of a de Bruijn digraph for all $k\ge 1$.

\begin{corollary}\label{cor2.3}
For $d\ge 2$, $\gamma_{k}(B(d,m))=\Big\lceil d^{m}\big/\sum_{j=0}^kd^j\Big\rceil$.
\end{corollary}
\noindent {\bf Proof.}
If $m\le k$, then clearly $\gamma_{k}(B(d,m))=\gamma_{k}(G_{B}(d^{m},d))=1=\big\lceil d^{m}\big/\sum_{j=0}^kd^j\big\rceil$ by Theorem \ref{thm2.2}, so the assertion holds. We may therefore assume
$m> k$. Let $m=ik+l$, where $i\geq 1$ and $0\leq l\le k-1$.
 Note that $d^m=(\sum_{j=0}^kd^j)(d^{m-k}-d^{m-k-1})+d^{m-k-1}$, $d^{m-k-1}=(\sum_{j=0}^kd^j)(d^{m-2k-1}-d^{m-2k-2})+d^{m-2k-2}$, $\cdots$,
then we have
\[ d^{m}=\left\{
\begin{array}{ll}
(\sum_{j=0}^kd^j)[(d^{m-k}-d^{m-k-1})+(d^{m-2k-1}-d^{m-2k-2})\\
~~~~~~+\cdots+(d^{m-(i-1)k-(i-2)}-d^{m-(i-1)k-(i-1)})]+d^{m-(i-1)k-(i-1)},\, \, \mbox{if $l< i$,}\\
(\sum_{j=0}^kd^j)[(d^{m-k}-d^{m-k-1})+(d^{m-2k-1}-d^{m-2k-2})\\
~~~~~~+\cdots+(d^{m-ik-(i-1)}-d^{m-ik-i})]+d^{m-ik-i}, \ \ \  \mbox{if $l\geq i$.}
\end{array}
\right.\]
Because $m=ik+l$ and $0\leq l\le k-1$, if $l< i$, then $d^{m-(i-1)k-(i-1)}=d^{l+k-(i-1)}\leq d^{k}$; and  if $l\geq i$, then $d^{m-ik-i}=d^{l-i}<d^{k}$. Thus
\[
\Big\lceil d^{m}\big/\sum_{j=0}^kd^j\Big\rceil=\left\{
\begin{array}{ll}
(d-1)(d^{m-k-1}+d^{m-2k-2}+\cdots+d^{m-(i-1)k-(i-1)})+1,\, \, \mbox{if $l<i$,}\\
(d-1)(d^{m-k-1}+d^{m-2k-2}+\cdots+d^{m-ik-i})+1, \,\, \mbox{if $l\geq i$.}
\end{array}
\right.\]
Hence  either $x=d^{m-k-1}+d^{m-2k-2}+\cdots+d^{m-(i-1)k-(i-1)}$ or $x=d^{m-k-1}+d^{m-2k-2}+\cdots+d^{m-ik-i}$ in $B(d,m)$ satisfies the congruence equation
$(d-1)x\equiv \big\lceil d^{m}\big/\sum_{j=0}^kd^j\big\rceil-h ~~({\rm mod} ~n)$
where  $h=1$ and $0\leq h\sum_{j=0}^{k-1}d^j\leq (\sum_{j=0}^kd^j)\lceil d^{m}\big/\sum_{j=0}^kd^j\rceil-d^{m}$.  Therefore, $\gamma_{k}(B(d,m))=\big\lceil d^{m}\big/\sum_{j=0}^kd^j\big\rceil$ by Theorem \ref{thm2.2}.
~$\Box$

As an application of Corollary \ref{cor2.2}, we provide a new sufficient condition for $\gamma_{k}(G_B(n,d))$ equal to $\big\lceil n\big/\sum_{j=0}^kd^j\big\rceil$. For this purpose,  we need the following result in elementary number theory.

For notational convenience,  $m\!\mid \!n$  means that $m$ divides
$n$ and $m\nmid n$  means that $m$ does not divide $n$ where $m, n$ are integers.  For integers $a_1, a_2, \ldots, a_n$,
the {\em greatest common divisor} of $a_1, a_2, \ldots, a_n$ is denoted by $(a_1, a_2, \ldots, a_n)$.

\begin{lemma}{\rm(\cite{pp})}\label{lem2.3}
For integers $a_1, a_2, \ldots, a_m$ $(m\ge 1)$, $b$ and  $n$, the congruence equation $\sum_{i=1}^ma_{i}x_{i}\equiv b$ $({\rm mod}~n)$ has at least a solution if
and only if $(a_{1}, a_{2}, \ldots, a_{m}, n)\!\mid \! b$.
\end{lemma}

\begin{theorem}\label{cor2.4}
For every generalized de Bruijn digraph $G_B(n,d)$, if both $n$ and $d$ satisfy one of the following conditions:

\noindent{\rm (i)} $\sum_{j=0}^kd^j\!\mid\! n$ and $(d-1,n)\!\mid\! n\big/\sum_{j=0}^kd^j$,

\noindent{\rm(ii)} $\big\lceil n\big/\sum_{j=0}^kd^j\big\rceil\equiv q$  {\rm(mod $(d-1, n)$)}, where $q$ satisfies
the inequality
$0\leq q(\sum_{j=0}^{k-1}d^j)\leq (\sum_{j=0}^kd^j)\big\lceil n\big/\sum_{j=0}^kd^j\big\rceil-n$,

\noindent then $\gamma_{k}(G_B(n,d))=\big\lceil n\big/\sum_{j=0}^kd^j\big\rceil$ and there is a vertex $x\in V(G_{B}(n,d))$ such that
$D=\{x, x+1, \cdots, x+\big\lceil n\big/\sum_{j=0}^kd^j\big\rceil-1\}$ is a consecutive minimum distance $k$-dominating set of $G_B(n,d)$.
\end{theorem}
\noindent {\bf Proof.}
Let $n$ and $d$ satisfy one of the conditions (i)-(ii). We show that
$G_B(n,d)$ contains a vertex $x$ such that $D=\{x, x+1, \cdots, x+\big\lceil n\big/\sum_{j=0}^kd^j\big\rceil-1\}$ is a consecutive minimum distance $k$-dominating set of $G_B(n,d)$. By Theorem \ref{thm2.2}, it suffices to show that there exists a vertex $x\in V(G_{B}(n,d))$ satisfies $(d-1)x\equiv \Big\lceil n\Big/\sum_{j=0}^kd^j\Big\rceil-h\, \,({\rm mod}\, n)$ (Eq. (\ref{formu6})) for some $h$ where $0\leq \big(\sum_{j=0}^{k-1}d^j\big)h\leq \big(\sum_{j=0}^kd^j\big)\lceil n\big/\sum_{j=0}^kd^j\rceil-n$.

(i) Suppose that $\sum_{j=0}^kd^j\!\mid\! n$ and $(d-1,n)\!\mid\! n\big/\sum_{j=0}^kd^j$. By Lemma \ref{lem2.3}, there is a vertex $x\in V(G_{B}(n,d))$ satisfying $(d-1)x\equiv n\big/\sum_{j=0}^kd^j\, ({\rm mod}\, n)$, so the assertion follows directly from Corollary \ref{cor2.2}.

(ii) Suppose that $\big\lceil n\big/\sum_{j=0}^kd^j\big\rceil\equiv q$  {\rm(mod $(d-1, n)$)}, where $q$ satisfies
the inequality
$0\leq q(\sum_{j=0}^{k-1}d^j)\leq (\sum_{j=0}^kd^j)\big\lceil n\big/\sum_{j=0}^kd^j\big\rceil-n$. Let $(d-1,n)=r$ and $\big\lceil n\big/\sum_{j=0}^kd^j\big\rceil=pr+q$ where $p\geq0$ and $0\le q\le r-1$. Set $q=h$.
Since $(d-1,n)| pr$,  the equation $(d-1)x\equiv pr \, \,({\rm mod}\, n)$ has a solution by Lemma \ref{lem2.3}. Hence,
there exists a vertex $x\in V(G_{B}(n,d))$ satisfying $(d-1)x\equiv \big\lceil n\big/\sum_{j=0}^kd^j\big\rceil-h \, \,({\rm mod}\, n)$, as desired. ~$\Box$

 By applying Theorems \ref{thm2.1} and \ref{thm2.2}, we obtain the following  sufficient condition for $\gamma_{k}(G_B(n,d))$ equal to $\big\lceil n\big/\sum_{j=0}^kd^j\big\rceil$.

\begin{theorem}\label{thm2.3}
 If $n=p(\sum_{j=0}^kd^j)+q$, where $p\geq 1$ and $1\leq q\leq {\rm min}\big\{1+2\sum_{j=0}^{k-1}d^j, \sum_{j=1}^{k}d^j\big\}$, then
$\gamma_{k}(G_{B}(n,d))=\big\lceil n\big/\sum_{j=0}^kd^j\big\rceil.$
\end{theorem}
\noindent {\bf Proof.}
By Theorem \ref{thm2.1}, we have known that
$G_B(n,d)$ contains a vertex  satisfying (1).
Let $x$ be such a vertex and let $D=\{x,x+1,\cdots,x+\big\lceil n\big/\sum_{j=0}^kd^j\big\rceil-1\}$.
We claim that $D$ is
a distance $k$-dominating set of $G_B(n,d)$. By the definition, it suffices to show that $\bigcup _{i=0}^{k}O_{i}(D)=V(G_B(n,d))$.

As before, we first show  the vertices of $O_{i-1}(D)\cup O_{i}(D)$ are consecutive for all $i, 1\leq i\le k$.
As already observed in Theorem \ref{thm2.2},   we have $O_{i}(D)=[d^{i}x, d^{i}\big(x+\big
\lceil n \big/\sum_{j=0}^kd^j\big\rceil\big)-1]\,\, ({\rm mod}\, n)$ and $|O_{i}(D)|=d^{i}\big
\lceil n \big/\sum_{j=0}^kd^j\big\rceil$ for all $i, 0\le i\le k$.
Since $x$ satisfies the inequality (\ref{formu1}), there exists an integer $h$, $0\le h\le d-2$
such that $dx= x+\big\lceil n\big/\sum_{j=0}^kd^j\big\rceil-h \,\,({\rm mod}\, n)$.
\begin{align*}
&d^{2}x= d\Big(x+\Big
\lceil n\Big/\sum_{j=0}^kd^j\Big\rceil\Big)-dh \,\, ({\rm mod}\, n),\\
&d^{3}x= d^{2}\Big(x+\Big
\lceil n\Big/\sum_{j=0}^kd^j\Big\rceil\Big)-d^{2}h \,\, ({\rm mod}\, n),\\
&~~~~~~~\vdots\\
&d^{k}x= d^{k-1}\Big(x+\Big
\lceil n\Big/\sum_{j=0}^kd^j\Big\rceil\Big)-d^{k-1}h\,\, ({\rm mod}\, n).
\end{align*}
Since $O_{i}(D)=[d^{i}x, d^{i}\big(x+\big
\lceil n \big/\sum_{j=0}^kd^j\big\rceil\big)-1]\,\, ({\rm mod}\, n)$ for all $i, 0\le i\le k$,
the vertices of
$O_{i-1}(D)\cap O_{i}(D)\neq \emptyset$ are consecutive for all $i, 1\leq i\le k$.

By the above fact,  we  show that $\bigcup _{i=1}^{k}O_{i}(D)$ contains all the vertices of $G_B(n,d)\setminus D$ by showing
the vertices of $O_k(D)\cup D$ are consecutive.
 We consider the last vertex in $O_k(D)$. Since $n=p(\sum_{j=0}^kd^j)+q$, $\big\lceil n\big/\sum_{j=0}^kd^j\big\rceil \sum_{j=0}^kd^j=n-q+\sum_{j=0}^kd^j$.
 Hence, by $dx= x+\big\lceil n\big/\sum_{j=0}^kd^j\big\rceil-h \,\,({\rm mod}\, n)$ where $0\leq h\leq d-2$, we have
\begin{align*}
d^{k}x+d^{k}\Big\lceil n\Big/\sum_{j=0}^kd^j\Big\rceil-1&= d^{k-1}\Big(x+\Big\lceil n\Big/\sum_{j=0}^kd^j\Big\rceil-h\Big)+d^{k}\Big\lceil n\Big/\sum_{j=0}^kd^j\Big\rceil-1\\
&=d^{k-1}x+(d^{k}+d^{k-1})\Big\lceil n\Big/\sum_{j=0}^kd^j\Big\rceil-d^{k-1}h-1\\
&=\cdots\\
%&= d^{k-2}\Big(x+\Big\lceil n\Big/\sum_{j=0}^kd^j\Big\rceil-h\Big)+(d^{k}+d^{k-1})
%\Big\lceil n\Big/\sum_{j=0}^kd^j\Big\rceil-d^{k-1}h-1\\
%&= d^{k-2}x+(d^{k}+d^{k-1}+d^{k-2})\Big\lceil n\Big/\sum_{j=0}^kd^j\Big\rceil-(d^{k-1}+d^{k-2})h-1\\
&= (x-1)+\Big\lceil n\Big/\sum_{j=0}^kd^j\Big\rceil \sum_{j=0}^kd^j-h\sum_{j=0}^{k-1}d^j\, \, ({\rm mod}\, n)\\
&= (x-1)+1+(d-h)\sum_{j=0}^{k-1}d^j-q\, \, ({\rm mod}\, n)\\
&\geq (x-1)+1+2\sum_{j=0}^{k-1}d^j-q \, \, ({\rm mod}\, n)\\
&\geq x-1,
\end{align*}
The last inequality holds, since $1\le q\le {\rm min}\big\{1+2\sum_{j=0}^{k-1}d^j, \sum_{j=1}^{k}d^j\big\}$.
Note that the vertices of $O_{i}(D)$ are  consecutive for all $i, 0\leq i\le k$, so $\bigcup _{i=1}^{k}O_{i}(D)\supseteq \{x+\big
\lceil n\big/\sum_{j=0}^kd^j\big\rceil,\ldots,  n-1,0,1,\ldots,  x-1\}$. This implies that $\bigcup _{i=1}^{k}O_{i}(D)\supseteq V(G_B(n,d))\setminus D$,
hence $D=\{x,x+1,x+2,\ldots, x+\big\lceil n\big/\sum_{j=0}^kd^j\big\rceil-1\}$ is a  distance $k$-dominating set of $G_B(n,d)$. Thus $\gamma_{k}(G_{B}(n,d))\le |D|=\big\lceil n\big/\sum_{j=0}^kd^j\big\rceil.$
By Theorem \ref{thm2.1}, $\gamma_{k}(G_{B}(n,d))=\big\lceil n\big/\sum_{j=0}^kd^j\big\rceil.$ ~$\Box$

\section{The minimum distance $k$-dominating sets in $G_K(n,d)$}
Tian and Xu \cite{tx} observed the following upper and lower bounds on $\gamma_{k}(G_K(n,d))$.

\begin{lemma}{\rm(\cite{tx})}\label{lem3.1}
 For any generalized Kautz digraph $G_K(n,d)$,
 $$\bigg\lceil n \Big/\sum_{j=0}^kd^j\bigg\rceil\le
\gamma_{k}(G_K(n,d))\le \bigg\lceil\frac{n}{d^{k}}\bigg\rceil.$$
\end{lemma}

In this section, we shall improve the above upper bound on $\gamma_{k}(G_K(n,d))$ by constructing a consecutive distance $k$-dominating set of $G_{K}(n,d)$.
%First we notice that the definitions of $G_{K}(n,d)$ and $i$-th out-neighborhood. For example, %$O(S)=\{15,14,13,12\}$, $O_{2}(S)=\{0,1, \cdots, 7\}$ and $O_{3}(S)=\{15,14,\cdots, 2\}$ for the subset $S=\{0,1\}$ %of $G_K(16,2)$.

\begin{theorem}\label{thm3.1}
 Let $G_K(n,d)$ be a generalized Kautz digraph. Then  $D=\big\{0,1,\cdots,\big\lceil n/(d^{k}+d^{k-1})\big\rceil-1\big\}$ is a distance $k$-dominating set of $G_{K}(n,d)$, and so
$$\gamma_{k}(G_K(n,d))\leq\bigg\lceil\frac{n}{d^{k}+d^{k-1}}\bigg\rceil.$$
\end{theorem}
\noindent {\bf Proof.}
 We show that $D$ is a distance $k$-dominating set of $G_{K}(n,d)$.
By the definitions of $G_{K}(n,d)$ and $i$-th out-neighborhood, if $k$ is odd, then we obtain
\begin{align*}
&O_{k-1}(D)=\big\{0,1,\cdots,d^{k-1}\big\lceil n/(d^{k}+d^{k-1})\big\rceil-1\big\},\\
 &O_{k}(D)=\big\{n-1,n-2,\cdots,n-d^{k}\big\lceil n/(d^{k}+d^{k-1})\big\rceil\big\};
\end{align*}
if $k$ is even, then
\begin{align*}
&O_{k-1}(D)=\big\{n-1,n-2,\cdots,n-d^{k-1}\big\lceil n/(d^{k}+d^{k-1})\big\rceil\big\},\\
&O_{k}(D)=\big\{0,1,\cdots,d^{k}\big\lceil n/(d^{k}+d^{k-1})\big\rceil-1\big\}.
\end{align*}
In both cases, we have $|O_{k-1}(D)|=d^{k-1}\big\lceil n/(d^{k}+d^{k-1})\big\rceil$ and
$|O_{k}(D)|=d^{k}\big\lceil n/(d^{k}+d^{k-1})\big\rceil$.
Note that the vertices of  $O_{k-1}(D)$ and  $O_{k}(D)$ are consecutive, and $(d^{k}+d^{k-1})\big\lceil n/(d^{k}+d^{k-1})\big\rceil\geq n$,  so $O_{k-1}(D)\cup O_{k}(D)=V(G_{K}(n,d))$. Hence  $D$ is a distance $k$-dominating set of $G_{K}(n,d)$. Therefore, $\gamma_{k}(G_K(n,d))\le |D|=
\big\lceil n/(d^{k}+d^{k-1})\big\rceil$.   ~$\Box$

\begin{remark}\label{rem3.1}{\rm
The upper bound on the distance $k$-domination number given in Theorem \ref{thm3.1} is sharp.
For example, we consider the digraph $G_K(7,2)$. We claim that $\gamma_{2}(G_K(7,2))=2=\big\lceil\frac{7}{2+4}\big\rceil$. Suppose not, we have $\gamma_{2}(G_K(7,2))=1$ by Lemma \ref{lem3.1}. Let $\{x_{0}\}$ be a minimum distance $2$-dominating set of $G_k(7,2)$.
Since $|O_i(x)|=d=2$ for each $x\in V(G_k(7,2))$, we have $O_{i}(x_{0})\cap O_{j}(x_{0})=\emptyset$ for all $0\leq i\neq j\leq 2$.
On the other hand, it can be verified that for each $x\in V(G_K(7,2))$, there exist integers $i,j$, $0\leq i\neq j\leq 2$, such that $O_{i}(x)\cap O_{j}(x)\neq\emptyset$ by the simply enumeration.
Thus each vertex $x$ of $G_K(7,2)$ can not form a distance $2$-dominating set of $G_K(7,2)$, as claimed.
By Theorem \ref{thm3.1},  $D=\{0,1\}$ must be a minimum distance $2$-dominating set of $G_K(7,2)$.}
\end{remark}

 The following result on the domination number of $G_K(n,d)$, due to Kikuchi and Shibata \cite{ks}, is an immediate consequence of Lemma \ref{lem3.1} and Theorem \ref{thm3.1}.

\begin{corollary}{\rm (\cite{ks})}\label{cor3.1}
 For every generalized Kautz digraph $G_K(n,d)$, $\gamma(G_K(n,d))=\big\lceil\frac{n}{d+1}\big\rceil$.
\end{corollary}

It seems to be  difficult to determine the minimum distance $k$-dominating
set for general generalized Kautz digraphs $G_{K}(n,d)$.  Now we present a sufficient condition for the distance $k$-domination number
of $G_{K}(n,d)$ to be the lower bound $\big\lceil n\big/\sum_{j=0}^kd^j\big\rceil$ in
Theorem \ref{thm3.1}.

\begin{theorem}\label{thm3.3}
For every generalized Kautz digraph $G_K(n,d)$, if $(d^{k-1}+d^{k})\big\lceil n\big/\sum_{j=0}^kd^j\big\rceil\geq n$ or $d^{k-1}\big\lceil n\big/\sum_{j=0}^kd^j\big\rceil\geq\big\lceil\frac{n}{d+1}\big\rceil$ then $\gamma_{k}(G_K(n,d))=\big\lceil n\big/\sum_{j=0}^kd^j\big\rceil$.
\end{theorem}
\noindent {\bf Proof.}
The proof is by directly constructing a (consecutive) distance $k$-dominating set of $G_K(n,d)$
with cardinality $\big\lceil n\big/\sum_{j=0}^kd^j\big\rceil$.
Let  $D=\big\{0,1,\cdots,\big\lceil n\big/\sum_{j=0}^kd^j\big\rceil-1\big\}$. We claim that $D$ is a distance $k$-dominating set of $G_{K}(n,d)$.
As  we have observed, if $k$ is odd, then
\begin{align*}
&O_{k-1}(D)=\Big\{0,1,\cdots,d^{k-1}\Big\lceil n\Big/\sum_{j=0}^kd^j\Big\rceil-1\Big\},\\ &O_{k}(D)=\Big\{n-1,n-2,\cdots,n-d^{k}\Big\lceil n\Big/\sum_{j=0}^kd^j\Big\rceil\Big\};
\end{align*}
 if $k$ is even,
then \begin{align*}
&O_{k-1}(D)=\Big\{n-1,n-2,\cdots,n-d^{k-1}\Big\lceil n\big/\sum_{j=0}^kd^j\Big\rceil\Big\},\\
&O_{k}(D)=\Big\{0,1,\cdots,d^{k}\Big\lceil n\Big/\sum_{j=0}^kd^j\big\rceil-1\Big\}.
\end{align*} Clearly,
$|O_{k-1}(D)|=d^{k-1}\big\lceil n\big/\sum_{j=0}^kd^j\big\rceil$ and
$|O_{k}(D)|=d^{k}\big\lceil n\big/\sum_{j=0}^kd^j\big\rceil$.

 Suppose that $(d^{k-1}+d^{k})\big\lceil n\big/\sum_{j=0}^kd^j\big\rceil\geq n$.  Note that the vertices of $O_{k-1}(D)$ and $O_{k}(D)$ are consecutive, so $O_{k-1}(D)\cup O_{k}(D)= V(G_K(n,d))$. Thus
$D=\big\{0,1,\cdots,\big\lceil n\big/\sum_{j=0}^kd^j\big\rceil-1\big\}$ is a distance $k$-dominating set of $G_{K}(n,d)$.

Suppose that $d^{k-1}\big\lceil n\big/\sum_{j=0}^kd^j\big\rceil\geq\big\lceil\frac{n}{d+1}\big\rceil$. By Lemma \ref{lem3.1} and Theorem \ref{thm3.1},  $D_{1}=\big\{0,1,\cdots,$ $\big\lceil\frac{n}{d+1}\big\rceil-1\big\}$ is a minimum dominating set of $G_K(n,d)$. Let $D_{1}'=\{n-1,n-2,\cdots,n-\big \lceil\frac{n}{d+1}\big\rceil \}$.
By the definition of $G_K(n,d)$, we have $O(D_{1}')=\{0,1,\cdots,d\big \lceil\frac{n}{d+1}\big\rceil-1 \}$.
Because $|D'_1\cup O(D'_1)|=(d+1)\big \lceil\frac{n}{d+1}\big\rceil\geq n$,
then $D_{1}'$ is also a minimum dominating set of $G_K(n,d)$.
 Since the vertices of $D$ are consecutive and $d^{k-1}\big\lceil n\big/\sum_{j=0}^kd^j\big\rceil\geq\big\lceil\frac{n}{d+1}\big\rceil$, we have either $O_{k-1}(D)\supseteq D_{1}$ or $O_{k-1}(D)\supseteq D_{1}'$. Hence $D=\big\{0,1,\cdots,\big\lceil n\big/\sum_{j=0}^kd^j\big\rceil-1\big\}$ is a distance $k$-dominating set of $G_{K}(n,d)$.
 ~$\Box$

\section{Closing remarks}

In this paper, we prove that the distance $k$-domination number of $G_B(n,d)$ takes on exactly  one of two values $\big\lceil n\big/\sum_{j=0}^kd^j\big\rceil$ and $\big\lceil n\big/\sum_{j=0}^kd^j\big\rceil+1$. In Theorems \ref{thm2.2}-\ref{thm2.3}, we provide  various sufficient conditions for $\gamma_{k}(G_B(n,d))$ equal to $\big\lceil n\big/\sum_{j=0}^kd^j\big\rceil$. It is of interest to determine the necessary and sufficient condition
for $\gamma_{k}(G_B(n,d))$ equal to $\big\lceil n\big/\sum_{j=0}^kd^j\big\rceil$. In Theorem \ref{thm3.1}, we establish the sharp upper bound on $\gamma_{k}(G_B(n,d))$. Furthermore, we provide a sufficient conditions for $\gamma_{k}(G_K(n,d))$ equal to $\big\lceil n\big/\sum_{j=0}^kd^j\big\rceil$ in Theorem \ref{thm3.3}.
We propose the following open problems.

\begin{prob}\label{prb1}
The sufficient condition in Theorem \ref{cor2.4} is also necessary for  $\gamma_{k}(G_B(n,d))$ equal to $\big\lceil n\big/\sum_{j=0}^kd^j\big\rceil$.
\end{prob}

For Problem \ref{prb1},  Dong, Shan and Kang \cite{dsk} proved that the assertion is true for the case when $k=1$.

\begin{prob}\label{prb2}
If $G_K(n,d)$ does not satisfy the conditions in Theorem \ref{thm3.3}, then  $\gamma_{k}(G_K(n,d))=\big\lceil n\setminus(d^{k-1}+d^{k})\big\rceil$.
\end{prob}

For Problem \ref{prb2}, if $k=1$, Corollary \ref{cor3.1},  due to Kikuchi and Shibata \cite{ks}, implies that the assertion is true.

\section*{ Acknowledgements}
The research was  supported in part by grants 11171207 and 11471210 of the
National Nature Science Foundation of China.
%Authors sincerely thanks the anonymous referees for their most valuable comment and suggestion.

\end{document}